\documentclass[11pt]{article}
\usepackage{amssymb}
\usepackage{graphicx}
\usepackage{amsfonts}
\usepackage{latexsym}
\usepackage{amsthm}
\usepackage{amsmath}

\usepackage{amssymb, amscd}

\usepackage{url}
\newtheorem{problem}{Problem}[section]
\newtheorem{theo}[problem]{Theorem}
\newtheorem{rem}[problem]{Remark}

\newtheorem{prop}[problem]{Proposition}
\newtheorem{cor}[problem]{Corollary}
\newtheorem{lema}[problem]{Lemma}

\begin{document}
\date{November 2004}
\title{{\huge \bf Equipartitions of measures in $\mathbb{R}^4$}}

\author{Rade  T. \v Zivaljevi\' c}%\footnote{The author was supported by the grant
%1643 of the Serbian Ministry of Science and Technology.}}

\maketitle

\begin{abstract}
A well known problem of B.~Gr\" unbaum \cite{Gru} asks whether for
every continuous mass distribution (measure) $d\mu = f\, dm$ on
$\mathbb{R}^n$ there exist $n$ hyperplanes dividing $\mathbb{R}^n$
into $2^n$ parts of equal measure. It is known that the answer is
positive in dimension $n=3$ \cite{Hadw} and negative for $n\geq
5$, \cite{Avis} \cite{Ramos}. We give a partial solution to Gr\"
unbaum's problem in the critical dimension $n=4$ by proving that
each measure $\mu$ in $\mathbb{R}^4$ admits an equipartition by
$4$ hyperplanes, provided that it is symmetric with respect to a
$2$-dimensional, affine subspace $L$ of $\mathbb{R}^4$. Moreover
we show, by computing the complete obstruction in the relevant
group of normal bordisms, that without the symmetry condition, a
naturally associated {\em topological problem} has a negative
solution. The computation is based on the Koschorke's exact
singularity sequence \cite{Kosch} and the remarkable properties of
the essentially unique, balanced, binary Gray code in dimension
$4$, \cite{Toot} \cite{Knu01}.
\end{abstract}

\section{Introduction}
A notorious open problem in geometric combinatorics and discrete
and computational geometry is the question whether each continuous
mass distribution $\mu$ in $\mathbb{R}^4$ admits an {\em
equipartition} by hyperplanes. This is an ``essentially
$4$-dimensional'' problem by the classification of V.~Klee
\cite{Klee}, indicating that the answer is known and positive in
all dimensions $\leq 3$ and negative in dimensions $\geq 5$.
Recall that a collection $H_1,H_2,\ldots,H_n$ of hyperplanes in
$\mathbb{R}^n$ is an equipartition for a mass distribution
(measure) $\mu$ if each of the $2^n$ ``orthants'' associated to
$\{H_j\}_{j=1}^n$ contains the fraction $1/2^n$ of the total mass.
In other words, $\mu(H^\epsilon)= 1/2^n \mu(\mathbb{R}^n)$ for
each $\epsilon =(\epsilon_1,\ldots, \epsilon_n)\in\{0,1\}^n$,
where $H^\epsilon:= H_1^{\epsilon_1}\cap\ldots\cap
H_n^{\epsilon_n}$ is the ``orthant'' associated to $\epsilon$ and
$H_i^0$ (respectively $H_i^{1}$) is the positive (respectively
negative) closed halfspace associated to the hyperplane $H_i$.

The progress in the general problem was slow after B.~Gr\" unbaum
posed the question in 1960, \cite{Gru}. H.~Hadwiger showed in
\cite{Hadw} that equipartitions exist for $n=3$. D.~Avis showed in
\cite{Avis} that there exist non-equipartitionable mass
distributions for $n\geq 5$. In the mean-time many related
questions were formulated and many of them solved \cite{BaMa2001}
\cite{BaMa2002} \cite{Mak2001} \cite{Ramos} \cite{VZ92}, a
connection with discrete and computational geometry was
established \cite{Yao-Yao} \cite{Yao} and the subject grew into a
separate branch of geometric combinatorics \cite{Ziv04}.
Nevertheless, the $4$-equipartition problem itself has resisted
all attempts and remains one of the central open problems in the
field.

In this paper we prove, Theorem~\ref{thm:plane}, that a measure
$\mu$ in $\mathbb{R}^4$ admits a $4$-equipartition if it is
symmetric with respect to a $2$-plane $L$ in $\mathbb{R}^4$.
Moreover we demonstrate in Theorem~\ref{thm:counterexample}, by
computing the obstruction in the relevant normal bordism group,
that there are no ``obvious'' topological obstacles for the
existence of non-equipartitionable measures.  This result may be
an indication that such a {\em peculiar} measure does exist  in
$\mathbb{R}^4$ and it is an intriguing question whether the
``topological counterexample'', provided by
Theorem~\ref{thm:counterexample}, can be turned into a genuine
counterexample.

\section{The CS/TM-scheme}

The {\em configuration space/test map} scheme \cite{Ziv04} has
emerged as one of the key principles for the application of
topological methods in geometric combinatorics and discrete and
computational geometry. The basic idea can be outlined as follows.

One starts with a configuration space or manifold $M_{\cal P}$ of
all candidates for the solution of a geometric/combinatorial
problem ${\cal P}$. The next step is a construction of a test map
$f : M_{\cal P}\rightarrow V_{\cal P}$ which measures how far is a
given candidate configuration $C\in M_{\cal P}$ from being a
solution. More precisely, there is a subspace $Z$ of the test
space $V_{\cal P}$ such that a configuration $C\in M_{\cal P}$ is
a solution if and only if $f(C)\in Z$. The inner symmetries of the
problem ${\cal P}$ typically show up at this stage. This means
that there is a group $G$ of symmetries of $X_{\cal P}$ which acts
on $V_{\cal P}$, such that $Z$ is a $G$-invariant subspace of
$V_{\cal P}$, which turns $f : M_{\cal P}\rightarrow V_{\cal P}$
into an equivariant map. If a configuration $C$ with the desired
property $f(C)\in Z$ does not exist, then there arises an
equivariant map $f : M_{\cal P}\rightarrow V_{\cal P}\setminus Z$.
The final step is to show by topological methods that such a map
does not exist.

The reader can follow the genesis of the method in review papers
\cite{Alon88} \cite{Bar93} \cite{Bjo91} \cite{Ziv2} \cite{Ziv3}
\cite{Ziv04} and see how the solutions of well known combinatorial
problems like Kneser's conjecture (L.~Lov\' asz \cite{Lov}), ``the
splitting necklace problem'' (N.~Alon \cite{Alon87}
 ), the Colored Tverberg problem (R.~\v Zivaljevi\' c,
S.~Vre\' cica, \cite{ZV2} \cite{VZ94}) etc.\ eventually led to the
formulation and codification of the general principle.

\subsection{The equipartition problem}
\label{sec:equi}

Our first choice for the configuration space suitable for the
equipartition problem is the manifold of all ordered collections
${\cal H} = (H_1,\ldots ,H_n)$ of oriented hyperplanes in
${\mathbb R}^n$.

Suppose that $e : \mathbb{R}^n\rightarrow \mathbb{R}^{n+1}$ is the
embedding defined by $e(x)= (x,1)$. Each oriented hyperplane $H$
in $e({\mathbb R}^n) \cong \{x\in {\mathbb R}^{n+1}\mid
x_{n+1}=1\}\subset \mathbb{R}^{n+1}$ is obtained as an
intersection $H = e({\mathbb R}^n)\cap H'$ for a unique  oriented,
$(n+1)$-dimensional subspace $H'\subset{\mathbb R}^{n+1}$. The
oriented subspace $H'$ is determined by the corresponding
orthogonal unit vector $u\in S^n\subset {\mathbb R}^{n+1}$, so the
natural environment for collections ${\cal H} = (H_1,\ldots,
H_n)$, and our second choice for the configuration space is
$M_{\cal P}:= (S^n)^n$. The group which acts on the configuration
manifold $M_{\cal P}$ is the reflection group $W_n :=
(\mathbb{Z}/2)^n\rtimes S_n$ where $S_n$ permutes the factors
while the subgroup $(\mathbb{Z}/2)^n$ is in charge of the
antipodal actions on individual spheres. The action of $W_n$ on
$M_{\cal P}:= (S^n)^n$ is not free, so our third choice for the
configuration space associated to the equipartition problem is
\begin{equation}
\label{eq:delta} M_{\cal P}^\delta = (S^n)^n_\delta :=\{x\in
(S^n)^n\mid x_i\neq\pm x_j \mbox{ {\rm for}  }  i\neq j\}.
\end{equation}
This space is a relative of the (standard) ``configuration space''
$F_m(S^n):=\{x\in (S^n)^m\mid x_i\neq x_j \mbox{ {\rm for}  }
i\neq j\}$ \cite{FaHu}. It has already appeared in Combinatorics,
for example in \cite{FeiZie}, where it is referred to as the
``signed configuration space''.

The associated ``orbit configuration space'' $(S^n)^n_\delta/W_n$
can be identified as a submanifold of the symmetric product
$SP^n(RP^n)$ of the projective space $RP^n$. This is the reason
why we occasionally denote this quotient by $SP_\delta(RP^n)$ and
view its elements as unordered collections of $n$ distinct lines
in $\mathbb{R}^{n+1}$.

The test space $V = V_{\cal P}$ is defined as follows.  If $\mu$
is a measure defined on ${\mathbb  R}^n$, let $\mu'$ be the
``push-down'' measure induced on ${\mathbb R}^{n+1}$ by the
embedding $e : {\mathbb  R}^n \hookrightarrow {\mathbb R}^{n+1}$,
 $\mu'(A):= \mu(e({\mathbb  R}^n)\cap A)$. A $n$-tuple $(u_1,\ldots
,u_n)\in (S^n)^n$ of unit vectors determines a $n$-tuple ${\cal H}
= (H_1',\ldots ,H_n')$ of oriented $(n+1)$-dimensional subspaces
of ${\mathbb  R}^{n+1}$. The $n$-tuple ${\cal H}$ dissects
${\mathbb R}^{n+1}$ into $2^n$-orthants $\mbox{\rm
Ort}_{\beta}({\cal H})$ which are naturally indexed by 0-1 vectors
$\beta \in \{0,1\}^n$. Let $b_{\beta} : M_{\cal P}\rightarrow
\mathbb{R}$ be the function defined by $b_{\beta}({\cal H}) :=
\mu'(\mbox{\rm Ort}_{\beta}({\cal H})) = \mu(\mbox{\rm
Ort}_{\beta}({\cal H})\cap e({\mathbb  R}^n))$. Let $B_{\mu} :
(S^n)^n\rightarrow {\mathbb R}^{2^n}$ be the function defined by
$B_{\mu}({\cal H}) = (b_{\beta}({\cal H}))_{\beta\in \{0,1\}^n}$.
The test space $V = V_{\cal P} \cong {\mathbb R}^{2^n}$ has a
natural action of the group $W_n := (\mathbb{Z}/2)^n\rtimes S_n$
such that the map $B_{\mu}$ is $W_n$-equivariant. Note that the
real $W_n$-representation $V$, restricted to the subgroup
$(\mathbb{Z}/2)^n \hookrightarrow W_n$, reduces to the regular
representation $\mbox{\rm Reg}((\mathbb{Z}/2)^n)$ of the group
$(\mathbb{Z}/2)^n$. The ``zero'' subspace $Z_{\cal P}$ is defined
as the trivial, $1$-dimensional $W_n$-representation $V_0$
contained in $V$. Let $U=U_n$ be the complementary
$W_n$-representation, $U_n \cong V/V_0$ and $A_{\mu} : (S^n)^n
\rightarrow U_n$ the induced, $W_n$-equivariant map. By the
construction we have the following proposition which says that
$A_{\mu}$ is a genuine test map for the $\mu$-equipartition
problem.
\begin{prop}
\label{prop:test} A $n$-tuple ${\cal H}=(H_1,\ldots ,H_n)\in
M_{\cal P}=(S^n)^n$ of oriented hyperplanes in $\mathbb{R}^n$ is
an equipartition of a measure $\mu$ defined on ${\mathbb  R}^n$ if
and only if $A_{\mu}({\cal H}) = 0$.
\end{prop}
\begin{cor}
\label{cor:equi} If there does not exist a $W_n$-equivariant map
$A : M_{\cal P}\rightarrow U_n\setminus\{0\}$, then each positive,
continuous mass distribution (measure) $d\mu = f\, dm$, where $dm$
is the Lebesgue measure, admits an equipartition by $n$
hyperplanes.
\end{cor}

\begin{rem} \label{rem:mere} {\rm The assumption
that $\mu$ is a measure absolutely continuous to the Lebesgue
measure on $\mathbb{R}^n$ is unnecessary restrictive. All we need
is the continuity of the test map $A_{\mu} : (S^n)^n \rightarrow
U_n$, a condition satisfied by a very large class of measures
having the desired continuity properties. Notable examples of
interesting measures that are not in this class are {\em counting}
measures $\nu_D$ of finite sets $D\subset \mathbb{R}^n$ defined by
$\nu_D(X):= \vert D\cap X\vert$. Note however, that all
equipartition results can be suitably extended to {\em weak
limits} of continuous measures, cf.\ \cite{MVZ}  for a general set
up.  An example of such a result applying to counting measures is
Corollary~\ref{cor:cloud} from Section~\ref{sec:plane}.}
\end{rem}

\subsection{Singular sets $\Sigma_\mu$}

It is a well known fact that for a free $G$-space $P$ and a
$G$-representation $V$, there does not exist a $G$-equivariant map
$f : P\rightarrow V\setminus\{0\}$ if and only if the associated
vector bundle $V\rightarrow P\times_GV\rightarrow P/G$ does not
admit a non-zero, continuous cross-section, cf.\ Proposition~I.7.2
in \cite{Dieck}.

A well known approach to the last question, applicable in the case
when $P$ is a free $G$-manifold, is the {\em singularity}\/
approach, cf.\ \cite{Kosch}. Given a $G$-map $h : P\rightarrow V$,
the {\em singularity set}\/ $\Sigma(h)$ of $h$ is the, possibly
empty, $G$-subspace of $P$ defined by $\Sigma(h):=h^{-1}(0)$. In
the case when $h$ is transverse to $0\in V$, the singularity set
$\Sigma(h)$ is a $G$-manifold.

If $h = A_\mu$ is the test map of a measure $\mu$, then the
associated singularity set $\Sigma_\mu := \Sigma(A_\mu)$ is
simply the set of all solutions to the equipartition problem for
$\mu$. In this case $\Sigma_\mu$ is often referred to as the {\em
solution set (manifold)} of $\mu$.

The singularity manifold $\Sigma(h)$ of a map $h$, sometimes
accompanied by the additional information recording the behavior
of $h$ in the tubular neighborhood of $\Sigma(h)$ (the normal
data), can be used for computation of an associated {\em
obstruction element}\/ in a suitable group of bordisms. The
singularity manifold and the associated normal data together yield
a very strong obstruction invariant which is in a number of
important cases complete in the sense that an equivariant map
exists if and only if these obstruction vanish. The reader is
referred to \cite{Kosch} for the general theory.

\subsection{Equipartitions of planar measures}
\label{sec:ravan}

The case $n=2$ of the equipartition problem is well known and
elementary. Nevertheless, we briefly review this case since it
serves as a fairly good illustration of general ideas in their
rudimentary form. According to the CS/TM-scheme, as presented in
Section~\ref{sec:equi}, the problem is to prove that there does
not exist a $W_2$-equivariant map $f : M_{\cal P}\rightarrow
U_2\setminus\{0\}$, where $W_2=\mathbb{D}_8$ is the dihedral
group, $M_{\cal P} = S^2\times S^2$, and $U_2$ the $3$-dimensional
real representation of $W_2$, described in Section~\ref{sec:equi}.

One can establish a slightly stronger statement that there does
not exist a $W_2$-equivariant map $f : M_{\cal P}^\delta
\rightarrow U_2\setminus\{0\}$ where $M_{\cal
P}^\delta=(S^2)^2_\delta = S^2\times S^2\setminus \{(x,y) \mid x=y
\mbox{ {\rm or} }x=-y\}$. The advantage of $(S^2)^2_\delta$ over
$(S^2)^2$ is that the former is a free $W_2$-space.

Let us see how the singularity approach works in the case of
planar equipartitions. For a {\em generic} measurable set
$A\subset\mathbb{R}^2$, the singularity $\Sigma_A$ of $A$, that is
the collection of all pairs $(L_1,L_2)$ of oriented lines in
$\mathbb{R}^2$ which form an equipartition for $A$, is a
$1$-dimensional $W_2$-manifold. For example if $A$ is a unit disc
$D$, the singularity $\Sigma_D$ is a union of $4$ circles. Here we
do not make precise what is meant by a generic measure. Instead we
naively assume, for the sake of this example, that there exists
such a notion of genericity for measurable sets/measures so that
each measurable set $A$ can be well approximated by generic
measures. Moreover, we assume that for any two measurable sets $A$
and $B$ there exists a path of generic measures $\mu_t, \, t\in
[0,1]$, so that $\mu_0$ is an approximation of $A$, $\mu_1$ is an
approximation for $B$ and the solution set
 \[
 \Sigma_{\{\mu_t\}_{t\in[0,1]}}:=\{(L_1,L_2;t)\mid (L_1,L_2) \mbox{
{\rm is an equipartition for } } \mu_t \} \subset (S^2)^2_\delta
\times [0,1]
 \]
is a $2$-dimensional manifold (bordism) connecting solution sets
for measures $\mu_0$ and $\mu_1$. The group $\Omega_1({\mathbb
D}_8)$ of classes of $1$-dimensional, free ${\mathbb
D}_8$-manifolds is isomorphic to $\mathbb{Z}/2\oplus \mathbb{Z}/2$
and the ${\mathbb D}_8$-solution manifold $\Sigma_D$, associated
to the unit disc in ${\mathbb  R}^2$, is easily shown to represent
a nontrivial element in this group.
\begin{figure}[hbt]
\centering
\includegraphics[scale=0.40]{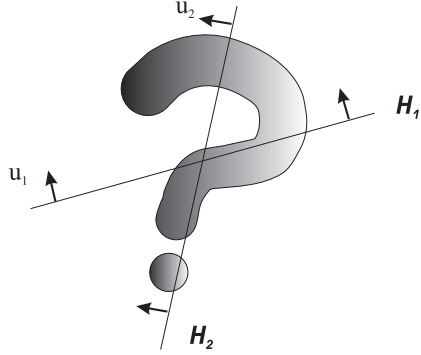}
\caption{An equipartition in the plane.} \label{fig:upitnik}
\end{figure}
It immediately follows that for any measurable set $A\subset
{\mathbb  R}^2$, or more generally for any continuous mass
distribution, the singularity $\Sigma_A$ is nonempty. Indeed,
suppose $\Sigma_A =\emptyset$. Let $\{\mu_t\}_{t\in[0,1]}$ be a
path of generic measures such that $\mu_0$ approximates $A$ and
$\mu_1$ approximates $D$. If these approximations are sufficiently
good, we deduce that the solution set $\Sigma_{\mu_0}$ is empty
and that $[\Sigma_{\mu_1}]$ and $[\Sigma_D]$ represent the same
element in $\Omega_1({\mathbb D}_8)$. This is a contradiction
since $\Sigma_{\mu_1} =
\partial(\Sigma)$ where $\Sigma = \Sigma_{\{\mu_t\}_{t\in[0,1]}}$, i.e.
$\Sigma_D$ would represent a trivial element in $\Omega_1({\mathbb
D}_8)$.

\bigskip
\begin{rem}
\label{rem1} {\rm It is worth noting that the scheme outlined
above, if applicable, shows that a general equipartition problem
can be solved by a careful analysis of the singularity set of a
well chosen, particular measure/measurable set, the unit disc
$D^2$ in our example above. Unfortunately, in higher dimensions
the unit balls do not represent generic measures, i.e. their
solution manifolds are very special and cannot be used for the
evaluation of the relevant obstruction elements. Instead, for this
purpose one can use measures distributed along a {\em convex
curve} $\Gamma$ in $\mathbb{R}^n$, Section~\ref{sec:convex2Gray}.
}

\end{rem}

\section{From convex curves to Gray codes}
\label{sec:convex2Gray}

\subsection{Convex curves}
\label{sec:convconv} A simple, smooth curve $\Gamma$ in
$\mathbb{R}^n$ is called {\em convex}\/ if the total multiplicity
of its intersection with an affine hyperplane does not exceed $n$.
Convex curves are classical objects appearing in many different
fields including the real algebraic geometry (``ecstatic points''
of curves), theory of convex polytopes (neighborly polytopes),
interpolations of functions (Tschebycheff systems), linear
ordinary differential equations (disconjugate equations) etc., see
\cite{Anis98} \cite{Arn96} \cite{Copp71} \cite{Scho54}
\cite{SeSh96} and the references in these papers. Standard
examples of convex curves are the ``moment curve'', or the
rational normal curve $M_n := \{(t,t^2,\ldots,t^n)\mid t\in
\mathbb{R}\}\subset \mathbb{R}^n$ and the standard trigonometric
curve $\Gamma_{2n}:=
\{(\cos{t},\sin{t},\cos{2t},\ldots,\cos{nt},\sin{nt})\mid t\in
[0,2\pi]\}\subset \mathbb{R}^{2n}$.

The importance of convex curves for the equipartition problem
stems from the fact that they minimize the number of intersections
with hyperplanes. As a consequence, each collection ${\cal
H}=\{H_1,\ldots, H_n\}$ of $n$ hyperplanes in $\mathbb{R}^n$ has
at most $n^2$ intersection points with $\Gamma$. It follows that
if $n=2d$ is even, and $\Gamma$ is a simple, closed convex curve,
then ${\cal H}$ divides $\Gamma$ in at most $n^2$ arcs. Suppose
that  $\mu$ is a measure concentrated on a closed, convex curve
$\Gamma\subset \mathbb{R}^{2d}$. Then if ${\cal H}$ is an
equipartition for $\mu$ then $n^2\geq 2^n$, i.e.\ $n\leq 4$. This
explains why the dimension $4$ is so special in this context.

\subsection{Gray codes}
\label{sec:GrayGray} 
\begin{figure}[hbt]
\centering
\includegraphics[scale=0.40]{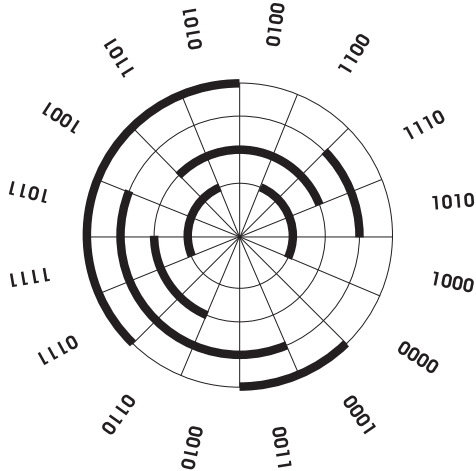}
\caption{The unique, balanced $4$-bit Gray code.} \label{slk:Gray}
\end{figure}

Gray codes arise in an attempt to describe
the solution manifold (singularity set) $\Sigma_\mu$ of a measure
$\mu$ distributed along a closed convex curve $\Gamma\subset
\mathbb{R}^4$. Our preferred example of such a curve is $\Gamma_4
= \{(z,z^2)\mid \vert z\vert =1\}\subset \mathbb{C}^2\cong
\mathbb{R}^4$. From here on we assume that $d\mu' = d\theta', \,
\theta' = {\rm arg}(z),$ is the ``arc length'' measure on the
circle $C = \{z\in \mathbb{C}\mid \vert z\vert =1\}$ and $d\mu=
d\theta$ the associated measure on $\Gamma_4$, where $\theta =
\theta'\circ\pi_1$ and $\pi_1: \mathbb{C}^2\rightarrow \mathbb{C}$
is the first projection. Similar analysis can be carried on for
other measures concentrated on $\Gamma_4$. Each collection $P =
\{p_1,p_2,p_3,p_4\}$ of $4$ points in $\Gamma_4$ is contained in a
unique hyperplane $H\subset \mathbb{R}^4$, hence the combinatorics
of $\mu$-equipartitions can be read off from the circle $C$, see
Figures~\ref{slk:Gray}, \ref{slk:jedankrug} and
\ref{slk:dvakruga}. As a consequence, an equipartition ${\cal
H}=\{H_1,H_2,H_3,H_3\}$ of the measure $d\theta$ is essentially a
division of the circle $C$ into $16$ arcs of equal length and
associating each of the $16$ division points $\{x_j\}_{j=0}^{15}$,
where $x_j=\epsilon^j~x_0$ and $\epsilon = \exp{2\pi i/16}$, to
one of hyperplanes $H_i$. Note that not all maps $\alpha :
\{x_j\}_{j=0}^{15}\rightarrow {\cal H}$ are allowed, cf.\
\cite{Ramos} p.\ 157. Each of the arcs $[x_j,x_{j+1}]$ belongs to
an orthant coded by an associated $4$-bit word
$\beta_j\in\{0,1\}^4$. Each of the $4$-bit words $\beta
\in\{0,1\}^4$ appears exactly once in the cyclic order inherited
from the order of intervals $[x_j,x_{j+1}]$. Moreover, moving from
one orthant to another along the curve $\Gamma_4$, that is from
the interval $[x_j,x_{j+1}]$ to the consecutive interval
$[x_{j+1},x_{j+2}]$, changes only one bit at a time. It follows
that sequence $\{\beta_j\}_{j=0}^{15}$ of $4$-bit words forms a so
called Gray code, \cite{Knu01} \cite{Ramos}. For a graph theorist,
a Gray code is a {\em Hamiltonian path} on a hypercube
$\{0,1\}^n$. For an engineer, a Gray code is a device useful for
converting digital signals into analog and vice versa,
\cite{Knu01}.

Our next observation is that Gray codes arising in the
equipartitions of the convex curve $\Gamma_4$ are quite special.
Indeed, each hyperplane $H_i$ intersects the convex curve
$\Gamma_4$ precisely $4$ times implying that the code must have
the same number of bit changes in each of four coordinate tracks.
Such codes are called {\em balanced}, Figure~\ref{slk:Gray}. A
Gray code with these properties was originally discovered by
G.C.~Tootill \cite{Toot}. One of its remarkable properties is that
it is unique up to permutation of coordinate tracks, see
\cite{Gilb} or \cite{Knu01}, Exercise~56. In particular if one
reads the code clockwise the same code is obtained, except that
the second and the third track interchange places, see
Figures~\ref{slk:Gray}--\ref{slk:dvakruga}.

There is one more attractive way to describe this code. As a
variation on the theme of the play ``Quad'' by S.~Beckett
\cite{Quad} where {\em ``... Four actors, whose colored hoods make
them identifiable yet anonymous, accomplish a relentless
closed-circuit drama ...''} (R.~Frieling), one can design a scheme
for a play based on the balanced Gray code. In this scheme the
stage begins and ends empty; $4$ actors enter and exit one at a
time, running through all $16$ possible subsets, and each actor is
supposed to enter (leave) the stage precisely $2$ times, see
\cite{Knu01} Exercise~65 (attributed to B.~Stevens) for a similar
idea.

\subsection{The solution set $\Sigma_{\theta}$}
\label{sec:solsol}

The analysis from Sections~\ref{sec:convconv} and
\ref{sec:GrayGray} allows us to describe the solution manifold
(singular set) $\Sigma_\theta$ of the measure $\theta$ on the
convex curve $\Gamma_4$, arising from the ``arc length'' measure
on the unit circle $C\subset \mathbb{C}$,
Section~\ref{sec:GrayGray}. This solution set is a $W_4$-invariant
subset of the configuration space $M_{\cal P}^\delta =
(S^4)^4_\delta$.

\begin{figure}[hbt]
 \centering
\includegraphics[scale=0.50]{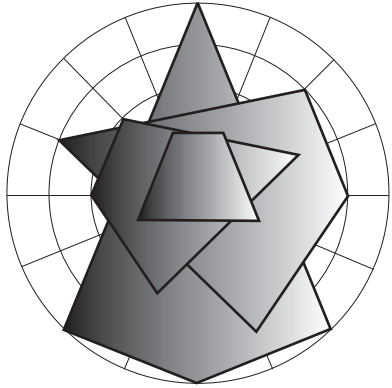}
\caption{An element of $\Sigma_{\theta}$...} \label{slk:jedankrug}
\end{figure}

Suppose that ${\cal H}=\{H_1,H_2,H_3,H_4\}\in\Sigma_\theta$. Each
hyperplane $H_j$ intersects the curve $\Gamma_4$ in four points,
vertices of a $3$-simplex $\sigma_j\subset H_j$. The image
$\pi_1(\sigma_j)$ of $\sigma_j$ by the projection map $\pi_1 :
\mathbb{C}^2\rightarrow \mathbb{C}$, sending the curve $\Gamma_4$
to the circle $C$, Section~\ref{sec:GrayGray}, is a convex polygon
in the plane. Figure~\ref{slk:jedankrug} displays these polygons
in the order corresponding to the chosen order of hyperplanes in
${\cal H}$.

By taking into account the (chosen) orientation of hyperplanes
$H_j$ and the induced orientations of simplices $\sigma_j\subset
H_j$, we arrive at the conclusion that Figure~\ref{slk:jedankrug}
is a fairly accurate description of an element ${\cal
H}=\{H_1,H_2,H_3,H_4\}\in\Sigma_\theta$. This element belongs to
an oriented circle $\Sigma\subset \Sigma_\theta$ of all solutions,
obtained essentially by rotating Figure~\ref{slk:jedankrug}
counterclockwise. All other circles in $\Sigma_\theta$ are
obtained by the action of the group $W_4$, in other words by
changing the orientations of simplices $\sigma_j$, performed by
the subgroup $(\mathbb{Z}/2)^4$, and by permuting the
circumcircles  of polygons in Figure~\ref{slk:jedankrug}. In other
words the permutation of polygons is achieved by permuting the
circles, the tracks in the balanced Gray code.

\begin{figure}[hbt]
 \centering
\includegraphics[scale=0.40]{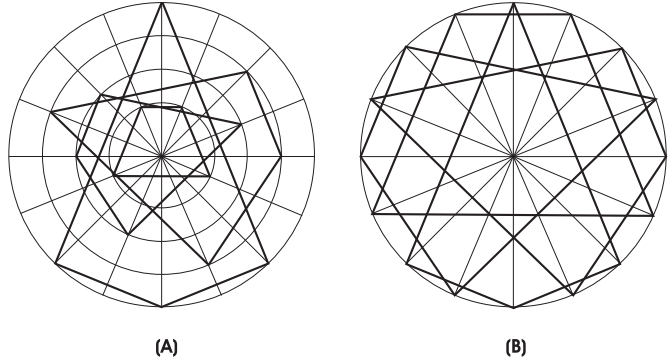}
\caption{... and the associated element in $\Sigma_{\theta}/W_4$.}
\label{slk:dvakruga}
\end{figure}

For example if we move Figure~\ref{slk:jedankrug} clockwise, we
obtain the circle $\xi(\Sigma)\subset\Sigma_\theta$ where $\xi :
[4]\rightarrow [4]$ is the permutation which keeps tracks $1$ and
$4$ fixed and interchanges tracks $2$ and $3$. We see this as a
manifestation of the symmetry of diagrams in
Figures~\ref{slk:jedankrug} and \ref{slk:dvakruga} with respect to
the vertical axes.

Let us turn to the question of non-degeneracy of the solution
manifold $\Sigma_\theta$. This is by definition the condition that
the associated test map $A_\theta : (S^4)^4_\delta \rightarrow
U_4$ is transverse to $0\in U_4$. Suppose that ${\mathcal
H}\in\Sigma_\theta$ and assume that $\{x_j\}_{j=0}^{15}$ are the
associated division points, Section~\ref{sec:GrayGray}. Then for a
small positive real number $\epsilon > 0$, the angles $y_j\in
(x_j-\epsilon, x_j+\epsilon)$ can be used as the coordinates on
$(S^4)^4_\delta$ in the neighborhood of ${\mathcal H}\in
(S^4)^4_\delta$. The associated tangent vectors are
$\partial/\partial y_j\in T_{\mathcal H}((S^4)^4_\delta)$.
Similarly, the functions $b_\beta$ introduced in
Section~\ref{sec:equi} are coordinates on $V$, consequently the
functions $c_j:= b_{\beta_{j+1}}-b_{\beta_{j}}$, where
$\{\beta_j\}_{j=0}^{15}$ is the sequence defined in
Section~\ref{sec:GrayGray}, are coordinates on $U_4$. The proof is
completed by the observation that $dA_\theta(\partial/\partial
y_j)=-\partial/\partial c_j$.

For the future reference we record an essential part of this
analysis in the following proposition.

\begin{prop}
\label{prop:sigma-theta} The solution manifold
$\Sigma_\theta\subset M_{\cal P}^\delta = (S^4)^4_\delta$ of all
equipartitions for the measure $d\theta$ on the convex curve
$\Gamma_4$ is a (non-degenerated) $1$-dimensional $W_4$-manifold
which has $\vert W_4\vert = 2^44!$ connected components. The
quotient manifold $\Sigma':=\Sigma_\theta/W_4$ is a circle in the
manifold $(S^4)^4_\delta/W_4\cong SP^4_\delta(RP^4)\subset
SP^4(RP^4)$, where $SP^m(X):= X^m/S_m$ is the symmetric product of
$X$ and $SP^m_\delta(X)$ its subspace of ``square-free'' divisors.
\end{prop}

\begin{rem}\label{rem:orbita}
{\rm Figure~\ref{slk:dvakruga} (B) symbolically represents an
element of the circle $\Sigma'=\Sigma_{\theta}/W_4$. The
orientations of polygons are forgotten and all $4$ circles,
representing different tracks in the balanced Gray code visible in
Figure~\ref{slk:dvakruga} (A), ``merged together'' in
Figure~\ref{slk:dvakruga} (B).   }
\end{rem}

\section{From Gray codes to normal bordisms}
\label{sec:Gray2bordisms}

According to Corollary~\ref{cor:equi}, the $4$-equipartition
problem is closely related to the question if there exists a
$W_4$-equivariant map $f : (S^4)^4_\delta \rightarrow
U_4\setminus\{0\}$. In turn, this is equivalent to the question if
the vector bundle
 \begin{equation}
 \label{eqn:bundle}
 {\cal E}: \quad U_4\rightarrow
(S^4)^4_\delta\times_{W_4}U_0\rightarrow (S^4)^4_\delta/W_4
 \end{equation}
admits a non-zero, continuous cross-section. A complete
topological obstruction $\omega$ for the existence of such a
section lives in the normal bordism group $\Omega_1(M;{\cal
E}-TM)$ \cite{Kosch}, where $TM$ is the tangent bundle of the
manifold $M := (S^4)^4_\delta/W_4$ and ${\cal E}$ is the bundle
(\ref{eqn:bundle}). This group can be computed from the
Koschorke's exact singularity sequence, \cite{Kosch} Theorem~9.3.
The computation of the obstruction $\omega\in\Omega_1(M;{\cal
E}-TM)$ is based on this sequence and the analysis of the
singularity $\Sigma_\theta$ of the measure $\theta$ distributed
along a closed, convex curve $\Gamma_4\subset \mathbb{R}^4$,
Section~\ref{sec:convex2Gray}.

\subsection{Koschorke's exact singularity sequence}
\label{sec:Kosch}

One of the consequences of Koschorke's exact singularity sequence
\cite[Section~7]{Kosch}, is a short exact singularity sequence,
\cite{Kosch} Theorem~9.3., involving low dimensional normal
bordism groups $\Omega_j(X;\phi)$ and
$\widetilde{\Omega}_j(X;\phi)$ for $j\leq 4$, where $\phi = \phi^+
- \phi^-$ is a virtual vector bundle over $X$. The final fragment
of this sequence has the following form,
\begin{equation}
\label{eqn:Koschorke} \stackrel{\delta_2}{\longrightarrow}
\Omega_2(X;\phi)\stackrel{f_2}{\longrightarrow}{\widetilde{\Omega}}_2(X;\phi)
\stackrel{\sigma\circ j_2}{\longrightarrow} \mathbb{Z}/2
\stackrel{\delta_1}{\longrightarrow} \Omega_1(X;\phi)
\stackrel{f_1}{\longrightarrow} \widetilde
{\Omega}_1(X;\phi)\longrightarrow 0 .
\end{equation}
We are interested in this sequence in the case when $X = M =
(S^4)^4_\delta/W_4$ and $\phi = \phi^+ - \phi^- = {\cal E}-TM$.
More precisely, our objective is to evaluate the obstruction
element $\omega\in\Omega_1(M;{\cal E}-TM)$ defined in
Section~\ref{sec:Gray2bordisms}.

\subsection{The image of $\delta_1$} \label{sec:image-delta}

According to \cite{Kosch}, the image $\delta_1(1)$ of the
generator $1\in \mathbb{Z}/2$ ``... can be represented by the unit
circle with constant map and the standard parallelization,
suitably stabilized ...''. In this section we show that the image
$\delta_1(1)$ coincides with the obstruction element
$\omega\in\Omega_1(M;{\cal E}-TM)$. This is a consequence of the
following proposition.

\begin{prop}\label{prop:delta}
$f_1(\omega)=0$.
\end{prop}

\medskip\noindent
{\bf Proof:} The configuration space $(S^4)^4_\delta$ is simply
connected, hence the projection $(S^4)^4_\delta\rightarrow
(S^4)^4_\delta/W_4$ is a universal covering map. Let $\Sigma$ be a
circle, connected component of the solution manifold
$\Sigma_\theta$, and let $\Sigma':= \Sigma_\theta/W_4\subset
(S^4)^4_\delta/W_4$. Suppose that $\xi$ is the orientation line
bundle of the virtual bundle ${\mathcal E}-TM$. Then, by
definition of $\widetilde{\Omega}_1(X;\phi)$, \cite{Kosch} \S~9,
the image $f_1(\omega)$ of $\omega$ in
$\widetilde{\Omega}_1(X;{\mathcal E}-TM)$ is determined by
$\Sigma'$ and the restriction $\xi\vert\Sigma'$ of $\xi$ on
$\Sigma'$.

Since the circle $\Sigma'$ can be lifted to the circle $\Sigma$,
we conclude that $\Sigma'$ is contractible, hence
$\xi\vert\Sigma'$ is a trivial bundle. It follows that there
exists a map $g: D^2\rightarrow M$ such that $\partial(D^2)$ is
mapped bijectively to $\Sigma'$ and extension of the bundle
$g^\ast(\xi\vert\Sigma')$ to trivial bundle over $D^2$. Hence,
$f_1(\omega)$ is a trivial element in
$\widetilde{\Omega}_1(X;{\mathcal E}-TM)$. \hfill$\square$

\subsection{The image of $\sigma\circ j_2$} \label{sec:image-sigma}

In this section we focus on the calculation of the image of the
map $\widetilde{\Omega}_2(X,\phi)\stackrel{\sigma\circ
j_2}{\longrightarrow} \mathbb{ Z}/2$ in the Koschorke's
singularity exact sequence (\ref{eqn:Koschorke}). By definition,
an element $\alpha = [N,g,or]\in \widetilde{\Omega}_2(X,\phi)$ is
mapped to $\sigma\circ\delta_2(\alpha) := g^{\ast}(w_2(\phi))[N]$,
where $w_2(\phi)$ is the second Stiefel-Whitney class of the
virtual bundle $\phi = \phi^{+}-\phi^{-}$ and $[N]\in H_2(N,
\mathbb{ Z}/2)$ is the fundamental class of the surface $N$. From
here we conclude that elements $\alpha = [N,g,or]$ such that both
$g^\ast(\phi^+)$ and $g^\ast(\phi^-)$ are trivial vector bundles
can be ignored, in particular we ignore those elements where $g$
is homotopic to a constant map. Recall the exact sequence
\begin{equation}\label{eq:unoriented}
\longrightarrow  {\cal N}_2\longrightarrow {\cal
N}_2(X)\longrightarrow H_2(X,\mathbb{ Z}/2)\longrightarrow 0
\end{equation}
 where ${\cal N}_2(X)$ is the group of unoriented bordisms and
${\cal N}_2:={\cal N}_2(\ast)$. It follows that in the evaluation
of the image of $\sigma\circ\delta_2$ we are allowed to pick a
representative $\alpha = [N,g,or]$ in each of the homology classes
$x\in H_2(X,\mathbb{ Z}/2)$.

 The following standard lemma reduces the calculation of
$w_2(\phi)[N]$ to the calculations of Stiefel numbers of
individual bundles $\phi^+$ and $\phi^-$.

\begin{lema}
\label{lema:stiefel}
  {\small
\[
\begin{array}{ccl}
w(\phi) & = & w(\phi^+)\cdot w(\phi^-)^{-1} \\
& = & (1+w_1(\phi^+)+w_2(\phi^+)+ \ldots
)(1+w_1(\phi^-)+w_2(\phi^-)+w_1(\phi^-)^2+ \ldots  ) \\
& = & 1+ A_1 + A_2 + \ldots , \mbox{ {\rm\it where} }
\end{array}
\]
}
\[
A_1 = w_1(\phi^+)+w_1(\phi^-) \mbox{ {\rm and} } A_2 =
w_2(\phi^+)+w_1(\phi^+)w_1(\phi^-)+w_1(\phi^-)^2 +w_2(\phi^-)
\]
are terms of graduation $1$ and $2$ respectively. Consequently if
$\phi = \phi^{+}-\phi^{-}$ is a virtual vector bundle over a
surface $N$, then
\[
w_2(\phi)[N] = w_2(\phi^+)[N] +
(w_1(\phi^+)w_1(\phi^-))[N]+w_1(\phi^-)^2[N] +w_2(\phi^-)[N].
\]
\end{lema}

In our case $X = SP^4_{\delta}(RP^4)$. The following lemma, an
easy consequence of Poincar\' e duality, allows us to search for
surfaces $N$ representing nontrivial homology $2$-classes in the
symmetric product $SP^4(RP^4)$.

\begin{lema} There is an isomorphism
$H_2(SP^4_{\delta}(RP^4))\longrightarrow H_2(SP^4(RP^4))$ of
homology groups, induced by the inclusion map
$SP^4_{\delta}(RP^4)\hookrightarrow SP^4(RP^4)$.
\end{lema}
It is easy to see that $H_2(SP^4(RP^4))\cong
H_2(SP^\infty(RP^4))$. By the Dold-Thom theorem, \cite{Hatch}
\[
SP^\infty (RP^4) \simeq K(\mathbb{ Z}/2,1)\times K(\mathbb{
Z}/2,2)\times K(\mathbb{ Z}/2,3)\times K(\mathbb{ Z}/2,4).
\]
From here we deduce that $H_2(SP^4(RP^4))\cong \mathbb{ Z}/2\oplus
\mathbb{ Z}/2$.

It is not difficult to describe surfaces $N_1$ and $N_2$, embedded
in $SP^4_{\delta}(RP^4)\subset SP^4(RP^4)$, representing the
generators of this homology group. Suppose that $e_1$ and $e'_1$
are two disjoint circles embedded in $RP^4$, both representing the
nontrivial element in $H_1(RP^4;\mathbb{ Z}/2)$. Let $e_2\cong
RP^2$ be a projective plane embedded in $RP^4$, representing the
generator of $H_2(RP^4;\mathbb{ Z}/2)\cong \mathbb{ Z}/2$. Finally
suppose that  $\ast_1, \ast_2, \ast_3$ are three distinct points
in $RP^4$ such that $\ast_i\notin e_1\cup e'_1\cup e_2$ for each
$i=1,2,3$. As usual, elements of the symmetric product $SP^m(X)$
are thought of as positive ``divisors'', i.e. commutative and
associative formal sums $D= n_1x_1+\ldots +n_kx_k$ where $n_i\in
\mathbb{N}, x_i\in X$ and $\sum_{i=1}^k = m$. By definition let
\begin{equation}\label{eq:surf}
N_1:= \ast_1 + \ast_2 + e_1 + e'_1   \quad {\rm and} \quad N_2 :=
\ast_1 + \ast_2 + \ast_3 + e_2.
\end{equation}
In other words $N_1\cong S^1\times S^1$ is a torus embedded in
$SP^4_\delta(RP^4)$ where $D\in N_1 \Leftrightarrow D = \ast_1 +
\ast_2 +x + y$ for some $x\in e_1$ and $y\in e'_1$. Similarly,
$N_2\cong RP^2$ consists of all divisors of the form $D = \ast_1 +
\ast_2 + \ast_3 + x$ for some $x\in e_2$.

\medskip
In our case $\phi^+ = {\mathcal E}$ and $\phi^- =
T(SP^4_\delta(RP^4))$. We focus our attention on the bundles
$\phi^+_i:= \phi^+\vert N_i$ and $\phi^-_i:= \phi^-\vert N_i$ for
$i=1,2$. Recall that ${\mathcal E}\cong (S^4)^4_\delta
\times_{W_4}U$ where $U=U_4$ is the $15$-dimensional
representation of $W_4$ described in Section~\ref{sec:equi}. If
$\pi : (S^4)^4_\delta\longrightarrow SP^4_\delta(RP^4)$ is the
projection map then $Z_i:= \pi^{-1}(N_i)$ is a free,
$W_4$-submanifold of $(S^4)^4_\delta$ and $\phi^+_i\cong
Z_i\times_{W_4}U$. It is not difficult to describe these
manifolds.

The connected component of $Z_1$ is a torus $T^2=S^1\times
S^1\subset (S^4)^4_\delta$ and the stabilizer of $T^2$ is the
group $H_1=\mathbb{Z}/2\oplus \mathbb{Z}/2\subset (\mathbb{
Z}/2)^{4}\subset W_4 = (\mathbb{Z}/2)^{4}\rtimes S_4$ where $H_1 =
\mathbb{ Z}/2\oplus \mathbb{ Z}/2$ acts on $T^2 = S^1\times S^1$
by the product action. It follows that $Z_1 \cong
T^2\times_{H_1}W_4$ and
 \[
\phi^+_1\cong Z_1\times_{W_4}U \cong
(T^2\times_{H_1}W_4)\times_{W_4}U \cong T^2\times_{H_1}U.
 \]
Similarly, the connected component of $Z_2$ is the sphere
$S^2\subset (S^4)^4_\delta$ and its stabilizer is the group
$H_2=\mathbb{Z}/2\subset (\mathbb{Z}/2)^{4}\subset W_4$ where
$H_2$ acts on $S^2$ by the antipodal action. It follows that
$Z_2\cong S^2\times_{H_2}W_4$ and
 \[
\phi^+_2\cong Z_2\times_{W_4}U \cong
(S^2\times_{H_2}W_4)\times_{W_4}U \cong S^2\times_{H_2}U.
 \]
Keeping in mind that the restriction ${\rm
Res}^{W_4}_{\mathbb{Z}^{\oplus 4}}(U)$ is the regular (real)
representation of $\mathbb{Z}^{\oplus 4}$, minus the trivial
$1$-dimensional representation, it is easy to identify the bundles
$\phi^+_1=T^2\times_{H_1}U$ and $\phi^+_2=S^2\times_{H_2}U$. As a
preliminary step, let us describe some canonical line bundles over
$N_1\cong T^2/H_1\cong S^1/\mathbb{ Z}/2\times S^1/\mathbb{
Z}/2\cong T^2$ and $N_2\cong S^2/H_2\cong RP^2$.

There are $4$ real, $1$-dimensional representations of $H_1 =
\mathbb{ Z}/2\oplus \mathbb{ Z}/2$. If $\omega_1$ and $\omega_2$
are the generators of $H_1$, then $L_{\epsilon_1\epsilon_2},
\epsilon_1,\epsilon_2\in\{-1,+1\}$, is the representation
characterized by the condition $\omega_i(v)=\epsilon_i v$ for each
$v\in L_{\epsilon_1\epsilon_2}$. Let
$\lambda_{\bar\epsilon_1\bar\epsilon_2}:=T^2\times_{H_1}L_{\epsilon_1\epsilon_2}$
be the associated line bundle where $\bar\epsilon_i\in\{0,1\}$ and
$(-1)^{\bar\epsilon_i} = \epsilon_i$. For example $\lambda_{00}$
is the trivial bundle, usually denoted by $\epsilon$, while
$\lambda_{11}$ is the Cartesian product of $2$ canonical bundles
over $RP^1\cong S^1$. Let $\gamma = \gamma_n$ be the canonical
line bundle over $RP^n$. The decompositions of bundles $\phi^+_1$
and $\phi^+_2$ into line bundles is recorded in the following
statement.

\begin{prop}
 \label{prop:record1}
 \[
\phi^+_1 = T^2\times_{H_1}U\cong \epsilon^{\oplus
3}\oplus\lambda_{01}^{\oplus 4}\oplus \lambda_{10}^{\oplus
4}\oplus \lambda_{11}^{\oplus 4}
 \]
 \[
\phi^+_2 =S^2\times_{H_2}U\cong \epsilon^{\oplus 7}\oplus
\gamma^{\oplus 8}.
 \]
\end{prop}

The next step is the identification of the restrictions
$\phi_1^{-} := TM\vert_{N_1}$ and  $\phi_2^{-} := TM\vert_{N_2}$
of the tangent bundle $TM = T(SP^4_\delta(RP^4))$ on the surfaces
$N_1$ and $N_2$ respectively.

\begin{lema}\label{lema:tangent}
For each point $D = p_1+p_2+p_3+p_4\in SP^4_\delta(RP^4)$ there is
an isomorphism
 \[
T_D(SP^4_\delta(RP^4))\cong \oplus_{i=1}^4 T_{p_i}(RP^4).
 \]
\end{lema}

\begin{lema}{\rm (\cite{M-S})}\label{lema:virtual}
Let $\gamma_n$ be the canonical line bundle over the real
projective space $RP^n$ and $\epsilon$ the trivial line bundle.
Then there is an equality of virtual vector bundles $T(RP^n)=
\gamma_n^{\oplus (n+1)} - \epsilon$. As a consequence, the total
Stiefel-Whitney class of the restriction bundle $\xi =
T(RP^n)\vert_{RP^m}$ is
 \[
w(\xi) = (1+t)^{n+1} = 1 + \binom{n+1}{1}t + \binom{n+1}{2}t^2 +
\ldots + \binom{n+1}{m}t^m.
 \]
\end{lema}

\begin{prop}\label{prop:record2}
\[
\begin{array}{cccll}
\phi^{-}_{1} & = & TM\vert_{N_1} & \cong & \epsilon^{\oplus 6}
\oplus \lambda_{01}^{\oplus
5} \oplus \lambda_{10}^{\oplus 5}\\
\phi^{-}_{2} & = & TM\vert_{N_2} & \cong & \epsilon^{\oplus 11}
\oplus \gamma^{\oplus 5}.
\end{array}
\]
\end{prop}

\medskip\noindent
{\bf Proof:} By Lemmas~\ref{lema:tangent} and \ref{lema:virtual},
there is an equality of virtual bundles

{\small\[
\begin{array}{cllll}
TM\vert_{N_1} & = & \epsilon^{\oplus 8} + (T(RP^4)\vert e_1 \times
T(RP^4)\vert e'_1) & = & \epsilon^{\oplus 8} + (\gamma_4^{\oplus
5}\vert {e_1} - \epsilon)\times
(\gamma_4^{\oplus 5}\vert {e'_1} - \epsilon)\\
& = & \epsilon^{\oplus 8} + \lambda_{01}^{\oplus 5} +
\lambda_{10}^{\oplus 5} - \epsilon^{\oplus 2} & = &
\epsilon^{\oplus 6} + \lambda_{01}^{\oplus 5} +
\lambda_{10}^{\oplus 5}.
\end{array}
\]}
Similarly,
\[
TM\vert_{N_2}\cong \epsilon^{\oplus 12} + T(RP^4)\vert e_2 =
\epsilon^{\oplus 12} + (\gamma_4^{\oplus 5}\vert e_2) - \epsilon)
= \epsilon^{\oplus 11} + \gamma^{\oplus 5}.
\]

\begin{cor}\label{cor:razlike}
\[
\begin{array}{cll}
\phi^{+}_{1} - \phi^{-}_{1} & = & \lambda_{11}^{\oplus 4} -
\lambda_{01} - \lambda_{10} - \epsilon^{\oplus 3}\\
\phi^{+}_{2} - \phi^{-}_{2} & = & \gamma^{\oplus 3} -
\epsilon^{\oplus 4}.
\end{array}
\]
\end{cor}

\begin{prop}
\label{cor:klase} ${\rm Im}(\sigma\circ j_2) = \mathbb{Z}/2$.
\end{prop}

\medskip\noindent
{\bf Proof:} It is a basic fact that
$H^\ast(RP^2;\mathbb{Z}/2)\cong (\mathbb{Z}/2)[t]/(t^3)$  and
$H^\ast(T^2;\mathbb{Z}/2)\cong\Lambda[a,b]$, where $\Lambda[a,b]$
is a $(\mathbb{Z}/2)$-exterior algebra generated by two elements
of degree $1$. The first two characteristic classes of line
bundles $\lambda_{\bar\epsilon_1\bar\epsilon_2}$ and $\gamma$,
defined over $N_1=T^2$ and $N_2=RP^2$ respectively, are
\begin{equation}\label{eqn:Stiefel}
w_1(\lambda_{\bar\epsilon_1\bar\epsilon_2}) = \bar\epsilon_1 a +
\bar\epsilon_2 b, \quad w_1(\gamma)=t, \quad
w_2(\lambda_{\bar\epsilon_1\bar\epsilon_2})= w_2(\gamma)=0.
\end{equation}
From here we deduce that the total Stiefel-Whitney classes of
bundles $\phi_1^+$ and $\phi_1^-$ are respectively
$w(\lambda_{11}^{\oplus 4})= 1$ and
$w(\lambda_{01}+\lambda_{10})=1+ a + b +ab$. This, together with
the Lemma~\ref{lema:stiefel} and the first equality from
Corollary~\ref{cor:razlike}, implies that
\begin{equation}\label{eqn:razlika}
w_2(\phi^{+}_{1} - \phi^{-}_{1})=w_1(\lambda_{01} +
\lambda_{10})^2 + w_2(\lambda_{01} + \lambda_{10})=ab.
\end{equation}
Similarly,
\begin{equation}\label{eqn:jos}
w_2(\phi^{+}_{2} - \phi^{-}_{2}) = w_2(\gamma^{\oplus 3})= t^2.
\end{equation}
In other words
\begin{equation}\label{eqn:slika}
w_2(\phi^{+}_{1} - \phi^{-}_{1})[N_1] = w_2(\phi^{+}_{2} -
\phi^{-}_{2})[N_2]=1
\end{equation}
and we finally conclude that ${\rm Im}(\sigma\circ j_2) =
\mathbb{Z}/2$. \hfill$\square$

\section{Results and proofs}
\subsection{Measures admitting a $2$-plane of symmetry}
\label{sec:plane}

In this section we show that each measure with a $2$-dimensional
plane of symmetry admits a $4$-equipartition.

\begin{theo}\label{thm:plane}
Suppose that $\mu$ is a measure on $\mathbb{R}^4$ admitting a
$2$-dimensional plane of symmetry in the sense that for some
$2$-plane $L\subset \mathbb{R}^4$ and the associated reflection
$R_L : \mathbb{R}^4\rightarrow \mathbb{R}^4$, for each measurable
set $A\subset \mathbb{R}^4$, $\mu(A) = \mu(R_L(A))$. Then $\mu$
admits a $4$-equipartition.
\end{theo}

\medskip\noindent
{\bf Proof:} Without loss of generality we assume that $L=
\mathbb{C}_{(2)}$ in the decomposition $\mathbb{R}^4\cong
\mathbb{C}^{2}\cong \mathbb{C}_{(1)}\oplus \mathbb{C}_{(2)}$. In
that case the reflection $R=R_L$ is the map described by
$R(z_1,z_2) = (-z_1,z_2)$, in particular $R$ is a symmetry of the
{\em convex curve} (Section~\ref{sec:convex2Gray}) $\Gamma_4
=\{(z,z^2)\in \mathbb{C}^2\mid \vert x\vert =1 \}$.

In pursuit of an equipartition of a general measure in
$\mathbb{R}^4$, we introduced in Section~\ref{sec:equi} the
configuration space $(S^4)^4_\delta = \{v\in S^4 \mid v_i\neq \pm
v_j \mbox{ {\rm for} } i\neq j \}$. Recall that the group $W_4 =
(\mathbb{Z}/2)^4\rtimes S_4$ acts freely on this configuration
space and that the problem of $4$-equipartitions was reduced,
Corollary~\ref{cor:equi}, to the question of the existence of a
$W_4$-equivariant map $f : (S^4)^4_\delta \rightarrow
U_4\setminus\{0\}$, where $U_4$ is the $15$-dimensional, real
representation of $W_4$ defined in Section~\ref{sec:equi}.

If the measure $\mu$ admits an additional symmetry, e.g.\ if it
admits a plane of symmetry $L$ such that the associated reflection
$R_L$ keeps $\mu$ invariant, it is natural to enlarge the group
$W_4$ by this transformation.

Assume as before that $\mathbb{R}^4$ is identified to the affine
subspace $\mathbb{R}^4 + e_5\subset \mathbb{R}^5$. The isometry
$R_L : \mathbb{R}^4\rightarrow \mathbb{R}^4$ is extended to the
unique isometry $\widehat{R}_L$ of $\mathbb{R}^5$ such that
$\widehat{R}_L(e_5)=e_5$.

Let $\mathbb{Z}/2$ be the group generated by $\widehat{R}_L$.
Define the enlarged group of symmetries as the direct product $G
:= W_4\times \mathbb{Z}/2$. The action of $W_4$ on
$(S^4)^4_\delta$ can be extended to the action of the group $G$ by
the requirement that $\widehat{R}_L(v_1,v_2,v_3,v_4) =
(\widehat{R}_L(v_1),\widehat{R}_L(v_2),\widehat{R}_L(v_3),\widehat{R}_L(v_4))$.
In order to make this action free, let us introduce an even
smaller configuration space
\begin{equation}\label{eq:small}
\begin{array} {ccl}
(S^4)^4_\Delta & := (S^4)^4_\delta \setminus {\cal F}
\end{array}
\end{equation}
where ${\cal F}\subset (S^4)^4_\Delta$ is the subset of points $v
= (v_1,v_2,v_3,v_4)$ such that the stabilizer ${\rm Stab}_G(v)$ is
a non-trivial group.

Note that $\widehat{R}_L$ is the reflection with respect to the
$3$-plane $L + \mathbb{R}e_5\subset \mathbb{R}^5$. Consequently,
$v\in (S^4)^4_\delta$ is in ${\cal F}$ if (up to a permutation of
coordinates $v_i$) either,
\[
(\widehat{R}_L(v_1) = v_2, \widehat{R}_L(v_3)=v_4) \quad {\rm
or}\quad (\widehat{R}_L(v_1) = v_1, \widehat{R}_L(v_2) = v_2,
\widehat{R}_L(v_3) = v_4).
\]

\medskip\noindent
{\bf Key observation:} The solution manifold $\Sigma_\theta\subset
(S^4)^4_\delta$ of all $4$-equipartitions of the convex curve
$\Gamma_4$, defined in Section~\ref{sec:solsol}, is
$\widehat{R}_L$-invariant. Moreover $\Sigma_\theta\subset
(S^4)^4_\Delta$, or in other words the $\mathbb{Z}/2$-action on
$\Sigma_\theta$ induced by $\widehat{R}_L$ is free. Indeed,
$\widehat{R}_L$ acts on the element displayed in
Figure~\ref{slk:jedankrug} by the rotation through the angle of
$180^\circ$.

\medskip\noindent
Denote by $\lambda$ both the non-trivial, $1$-dimensional real
representation of $\mathbb{ Z}/2$ and the associated,
$1$-dimensional $G$-representation induced by the projection
homomorphism $G\rightarrow \mathbb{ Z}/2$. Similarly, $U_4$ is
both the $15$-dimensional real $W_4$-representation defined in
Section~\ref{sec:equi} and the representation induced by the
projection homomorphism $G\rightarrow W_4$.

According to results of Section~\ref{sec:equi}, for the proof of
the theorem it is sufficient to show that a $W_4$-equivariant map
$f : (S^4)^4_\delta\rightarrow U_4$ must have a zero. This is a
consequence of the following stronger result.

\medskip\noindent
{\bf Claim:} There does not exist a $G$-equivariant map $f:
(S^4)^4_\Delta \rightarrow S(U_4\oplus\lambda)$, where
$S(U_4\oplus\lambda)$ is the $G$-invariant unit sphere in
$U_4\oplus\lambda$. In other words each $G$-invariant map $f:
(S^4)^4_\Delta \rightarrow  U_4\oplus\lambda$ has a zero.

\medskip\noindent
{\bf Proof of the Claim:} The claim is equivalent to the statement
that the vector bundle $\xi : (S^4)^4_\Delta\times_G
(U_4\oplus\lambda) \rightarrow  (S^4)^4_\Delta/G$ does not admit a
non-zero continuous cross-section. For this it is sufficient to
show that the top Stiefel-Whitney class $w_n(\xi)$ is non-zero. By
duality this is equivalent to the fact that for some (each) smooth
cross-section $s : (S^4)^4_\Delta/G\rightarrow
(S^4)^4_\Delta\times_G (U_4\oplus\lambda)$, transverse to the zero
section $Z$, the number of elements in $s^{-1}(Z)\subset
(S^4)^4_\Delta/G$ is odd. In the language of equivariant maps,
this is equivalent to the statement that for some smooth,
$G$-equivariant map $s =(\phi,\psi): (S^4)^4_\Delta \rightarrow
U_4\oplus\lambda$, transverse to $0\in E\oplus\lambda$, the number
of $G$-orbits in $s^{-1}(0)$ is odd.

Define $\phi : (S^4)^4_\Delta \rightarrow U_4$ as the restriction
of the test map $A_\theta : (S^4)^4_\delta \rightarrow U_4$
introduced in Section~\ref{sec:equi}, where $d\theta$ is the ``arc
length''-measure on $\Gamma_4 = \{(z,z^2)\in \mathbb{C}^2\mid
\vert z\vert =1\}$ introduced in Section~\ref{sec:GrayGray}.

Since $A_\theta$ is $W_4$-equivariant and $A_\theta\circ
\widehat{R}_L = A_\theta$, we conclude that $\phi$ is
$G$-equivariant. The orbit space $SP^4_\Delta(RP^4):=
(S^4)^4_\Delta/W_4$ is a smooth manifold. Note,
Section~\ref{sec:convex2Gray}, that the unique balanced Gray code
allowed us to identify the ``solution manifold'' of all
equipartitions of the curve $\Gamma_4\subset \mathbb{R}^4$ as the
circle $\Sigma'\subset SP^4_\Delta(RP^4)$. Moreover, this solution
manifold was shown to be non-degenerated,
Proposition~\ref{prop:sigma-theta}, in the sense that the
associated test map $A_\theta$ is transverse to $0\in U_4$.

The $\mathbb{ Z}/2$-action on $SP^4_\Delta(RP^4)$, induced by the
involution $\widehat{R}_L$, is free. Let $\psi' :
SP^4_\Delta(RP^4) \rightarrow \lambda$ be a $\mathbb{
Z}/2$-equivariant, smooth map such that $0\in\lambda$ is not a
critical value of the restriction $\psi'':=\psi'\vert \Sigma'$ of
$\psi'$ on the circle $\Sigma'$. Define $\psi$ as the composition
of $\psi'$ with the natural projection map
$(S^4)^4_\Delta\rightarrow SP^4_\Delta(RP^4)$.

Then it is not difficult to check that $s=(\phi,\psi)$ is a
smooth, $G$-equivariant map such that $0\in\lambda$ is not one of
its critical values.

Let us show that the number of $G$-orbits in the set $s^{-1}(0)$
is always an odd number. Note that $s^{-1}(0) = \phi^{-1}(0)\cap
Z(\psi)$ where $Z(\psi):=\psi^{-1}(0)$ is the zero set of $\psi$.
Since $\phi^{-1}(0)/G = \phi^{-1}(0)/W_4 = \Sigma'$, we observe
that the number of $G$-orbits in the set $s^{-1}(0)$ is equal to
the number of $\mathbb{ Z}/2$-orbits in the zero set $Z(\psi'')$
of the $\mathbb{Z}/2$-equivariant map $\psi'' : \Sigma'\rightarrow
\lambda$. Since $0$ is not a critical value of $\psi''$, the proof
is completed by an elementary observation that a $\mathbb{
Z}/2$-equivariant map $p : S^1\rightarrow \lambda$, transverse to
$0\in \lambda$, must have an odd number of $\mathbb{Z}/2$-orbits.

\hfill $\square$

\begin{cor}\label{cor:cloud}
Suppose that $D$ is a finite set of $16d$ distinct points  in
$\mathbb{R}^4$ which is symmetric with respect to a $2$-plane
$L\subset\mathbb{R}^4$. Then there exists a collection ${\cal H} =
\{H_1,H_2,H_3,H_4\}$ of distinct hyperplanes such that each of
$16$ associated open orthants contains not more than $d$ points
from the set $D$. Moreover, if $D$ is in general position in the
sense that no $5$ points belong to the same hyperplane, than each
of the open orthants contains at least $d-16$ elements from $D$.
In particular, if $D$ is a not necessarily symmetric set of $16d$
distinct points in $\mathbb{R}^4$, then for some collection ${\cal
H}$ of hyperplanes each of the associated open orthants contains
not more than $2d$ points from $D$.
\end{cor}

\begin{rem}
\label{rem:accident} {\rm Note that the proof of
Theorem~\ref{thm:plane} is based on the ``accident'' that the
convex curve $\Gamma_4$ and the associated solution manifold
$\Sigma'$ are $\mathbb{Z}/2$-spaces, where $\mathbb{Z}/2$ is the
group generated by the reflection $R_L$. As a consequence
$\Sigma'$ determines a nontrivial element in the group
$\Omega(\mathbb{Z}/2)$ of $\mathbb{Z}/2$-bordisms \cite{CoFl}
which means that one should be able to repeat the pattern of the
proof of the equipartition theorem in the planar case,
Section~\ref{sec:ravan}. Our proof essentially follows this idea.
For example the definition of $\psi'$ corresponds to the choice of
a non-trivial $\mathbb{Z}/2$-equivariant, line bundle on the
circle $\Sigma'$, a step used in the proof that $\Sigma'$ defines
a non-trivial element in $\Omega(\mathbb{Z}/2)$.}
\end{rem}

\subsection{Measures admitting a center or a $3$-plane of symmetry}
\label{sec:center}

For completeness we include proofs of $4$-equipartition results
for measures in $\mathbb{R}^4$ which admit a center or a $3$-plane
of symmetry, Proposition~\ref{prop:center} and
Proposition~\ref{prop:3plane}. These results formally resemble
Theorem~\ref{thm:plane} but there are important differences. Both
Propositions~\ref{prop:center} and \ref{prop:3plane} are easily
deduced from Hadwiger's result about $3$-equipartitions of
measures in $\mathbb{R}^3$, Theorem~\ref{thm:Hadwiger}. Contrary
to this, the proof of Theorem~\ref{thm:plane} relies on the
existence of a unique, ballanced, $4$-bit binary Gray code, so
this result appears to be an essential $4$-dimensional phenomenon.

\begin{theo}{\rm (\cite{Hadw})}\label{thm:Hadwiger}
Each continuous mass distribution $\mu$ in $\mathbb{R}^3$ admits a
$3$-equipartition, that is a collection $H_1, H_2, H_3$ of three
planes in $\mathbb{R}^3$ dissecting the ambient space into eight
octants of equal measure. Moreover, the first of these planes can
be chosen to contain arbitrary two points $A,B\in \mathbb{R}^3$
prescribed in advance.
\end{theo}

\begin{rem}
{\rm Hadwiger \cite{Hadw} deduced Theorem~\ref{thm:Hadwiger} from
the result that any two measures in $\mathbb{R}^3$ admit a
simultaneous equipartition by $2$ hyperplanes. Both results of
Hadwiger can be proved, along the lines of the CS/TM-scheme, by an
analysis of the set of all equipartitions for measures with
compact support, concentrated on the convex curve
$M_3=\{(t,t^2,t^3)\mid t\in \mathbb{R}\}$, see \cite{Ziv3}
Proposition~4.9. }
\end{rem}

\begin{prop}{\rm (\cite{Zie})}
\label{prop:center} Suppose that $\mu$ is a continuous mass
distribution in $\mathbb{R}^4$ which has a center of symmetry $O$.
Then $\mu$ admits a $4$-equipartition. Moreover, one of the
hyperplanes can be choosen in advance as  an arbitrary $3$-plane
passing through the center of symmetry $O$.
\end{prop}
\medskip\noindent
{\bf Proof:} Choose a ``halving'' hyperplane $H_1$ for $\mu$. One
can assume that $O\in H_1$. Let $\pi : \mathbb{R}^4\rightarrow
H_1$ be the orthogonal projection and let $\nu$ be the measure on
$H_1$ defined by $\nu(A) := \mu(\pi^{-1}(A)\cap H_1^+)$ where
$H_1^+$ is a closed halfspace bounded by $H_1$. Find $2$-planes
$P_2,P_3,P_4$ in $H_1$ which form a $3$-equipartition for $\nu$
such that $O\in P_i$ for each $i$. This is always possible by
Theorem~\ref{thm:Hadwiger}. Then ${\cal H}=\{H_1,H_2,H_3,H_4\}$ is
a $4$-equipartition for $\mu$ where $H_j:=\pi^{-1}(P_j)$ for
$j\geq 2$ is the hyperplane orthogonal to $H_1$ at $P_j$.
\hfill$\square$

\begin{cor}
There does not exist a centrally symmetric, convex closed curve in
$\mathbb{R}^4$.
\end{cor}
\medskip\noindent
{\bf Proof:} The ``arc length''-measure on such a curve would be
centrally symmetric. According to Proposition~\ref{prop:center},
its space of $4$-equipartitions is ``fibered'' over the
Grassmannian of all affine $3$-planes in $\mathbb{R}^4$, hence it
is at least $4$-dimensional. This is a contradiction. Indeed, by
the results of Section~\ref{sec:convex2Gray}, the solution
manifold of all $4$-equipartitions of a measure concentrated on a
closed, convex curve in $\mathbb{R}^4$ must be $1$-dimensional.
\hfill$\square$

\begin{prop}
\label{prop:3plane} Suppose that $\mu$ is a continuous mass
distribution in $\mathbb{R}^4$ which has a $3$-plane of symmetry.
Then there exists a $4$-equipartition for $\mu$.
\end{prop}
\medskip\noindent
{\bf Proof:} The proof is similar to the proof of
Proposition~\ref{prop:center}, so we omit the details.
\hfill$\square$

\subsection{Topological ``counterexample''}

\begin{theo}
\label{thm:counterexample} There exists a $W_4$-equivariant map $f
: (S^4)^4_\delta\rightarrow S(U_4)$.
\end{theo}

\medskip\noindent
{\bf Proof:} By \cite{Kosch} $\S 3$, a $W_4$-equivariant map $f :
(S^4)^4_\delta\rightarrow S(U_4)$ exists if and only if the
obstruction $\omega\in\Omega_1(M;{\cal E}-TM)$ in the
corresponding group of normal bordisms vanishes, cf.\
Section~\ref{sec:Gray2bordisms}. By Proposition~\ref{prop:delta},
$\omega$ is in the image of the map $\delta_1$. By
Proposition~\ref{cor:klase} the map $\sigma\circ j_2$ is onto and
the result follows from the exactness of the sequence
(\ref{eqn:Koschorke}). \hfill$\square$

\end{document}